\documentclass[12pt]{amsart}
\usepackage{amsfonts,amssymb,color}
\usepackage[mathscr]{eucal}
\usepackage{amsmath, amsthm}
\usepackage{mathrsfs}
\usepackage{amsbsy}
\usepackage{dsfont}
\usepackage{bbm}
\usepackage{wasysym}
\usepackage{stmaryrd}
\input xypic
\xyoption{all}
\begin{document}
\baselineskip = 5mm
\newcommand \ZZ {{\mathbb Z}} 
\newcommand \FF {{\mathbb F}} %
\newcommand \NN {{\mathbb N}} 
\newcommand \QQ {{\mathbb Q}} 
\newcommand \RR {{\mathbb R}} 
\newcommand \CC {{\mathbb C}} 
\newcommand \PR {{\mathbb P}} 
\newcommand \AF {{\mathbb A}} 
\newcommand \uno {{\mathbbm 1}}
\newcommand \Le {{\mathbbm L}}
\newcommand \bcA {{\mathscr A}}
\newcommand \bcB {{\mathscr B}}
\newcommand \bcC {{\mathscr C}}
\newcommand \bcD {{\mathscr D}}
\newcommand \bcE {{\mathscr E}}
\newcommand \bcF {{\mathscr F}}
\newcommand \bcG {{\mathscr G}}
\newcommand \bcH {{\mathscr H}}
\newcommand \bcM {{\mathscr M}}
\newcommand \bcN {{\mathscr N}}
\newcommand \bcI {{\mathscr I}}
\newcommand \bcK {{\mathscr K}}
\newcommand \bcL {{\mathscr L}}
\newcommand \bcO {{\mathscr O}}
\newcommand \bcP {{\mathscr P}}
\newcommand \bcQ {{\mathscr Q}}
\newcommand \bcR {{\mathscr R}}
\newcommand \bcS {{\mathscr S}}
\newcommand \bcT {{\mathscr T}}
\newcommand \bcU {{\mathscr U}}
\newcommand \bcV {{\mathscr V}}
\newcommand \bcW {{\mathscr W}}
\newcommand \bcX {{\mathscr X}}
\newcommand \bcY {{\mathscr Y}}
\newcommand \bcZ {{\mathscr Z}}
\newcommand \Spec {{\rm {Spec}}}
\newcommand \Pic {{\rm {Pic}}}
\newcommand \Jac {{{J}}}
\newcommand \Alb {{\rm {Alb}}}
\newcommand \NS {{{NS}}}
\newcommand \Corr {{Corr}}
\newcommand \Sym {{\rm {Sym}}}
\newcommand \Alt {{\rm {Alt}}}
\newcommand \Prym {{\rm {Prym}}}
\newcommand \cone {{\rm {cone}}}
\newcommand \cha {{\rm {char}}}
\newcommand \tr {{\rm {tr}}}
\newcommand \alg {{\rm {alg}}}
\newcommand \im {{\rm im}}
\newcommand \Hom {{\rm Hom}}
\newcommand \colim {{{\rm colim}\, }} 
\newcommand \End {{\rm {End}}}
\newcommand \coker {{\rm {coker}}}
\newcommand \id {{\rm {id}}}
\newcommand \tor {{\rm {tor}}}
\newcommand \spc {{\rm {sp}}}
\newcommand \Ob {{\rm Ob}}
\newcommand \Aut {{\rm Aut}}
\newcommand \cor {{\rm {cor}}}
\newcommand \res {{\rm {res}}}
\newcommand \Gal {{\rm {Gal}}}
\newcommand \PGL {{\rm {PGL}}}
\newcommand \Gr {{\rm {Gr}}}
\newcommand \Bl {{\rm {Bl}}}
\newcommand \Sing {{\rm {Sing}}}
\newcommand \spn {{\rm {span}}}
\newcommand \Nm {{\rm {Nm}}}
\newcommand \inv {{\rm {inv}}}
\newcommand \codim {{\rm {codim}}}
\newcommand \ptr {{\pi _2^{\rm tr}}}
\newcommand \gom {{\mathfrak m}}
\newcommand \goT {{\mathfrak T}}
\newcommand \goC {{\mathfrak C}}
\newcommand \goD {{\mathfrak D}}
\newcommand \goM {{\mathfrak M}}
\newcommand \goS {{\mathfrak S}}
\newcommand \goH {{\mathfrak H}}
\newcommand \sg {{\Sigma }}
\newcommand \CHM {{\mathscr C\! \mathscr M}}
\newcommand \DM {{\sf DM}}
\newcommand \FS {{FS}}
\newcommand \MM {{\mathscr M\! \mathscr M}}
\newcommand \HS {{\mathscr H\! \mathscr S}}
\newcommand \MHS {{\mathscr M\! \mathscr H\! \mathscr S}}
\newcommand \Vect {{\mathscr V\! ect}} %
\newcommand \Gm {{{\mathbb G}_{\rm m}}}
\newcommand \trdeg {{\rm {tr.deg}}}
\newcommand \znak {{\natural }}
\newtheorem{theorem}{Theorem}
\newtheorem{lemma}[theorem]{Lemma}
\newtheorem{corollary}[theorem]{Corollary}
\newtheorem{proposition}[theorem]{Proposition}
\newtheorem{remark}[theorem]{Remark}
\newtheorem{definition}[theorem]{Definition}
\newtheorem{conjecture}[theorem]{Conjecture}
\newtheorem{example}[theorem]{Example}
\newtheorem{question}[theorem]{Question}
\newtheorem{assumption}[theorem]{Assumption}
\newtheorem{fact}[theorem]{Fact}
\newtheorem{crucialquestion}[theorem]{Crucial Question}
\newcommand \lra {\longrightarrow}
\newcommand \hra {\hookrightarrow}
\def\blue {\color{blue}}
\def\red {\color{red}}
\def\green {\color{green}}
\newenvironment{pf}{\par\noindent{\em Proof}.}{\hfill\framebox(6,6)
\par\medskip}
\title[Algebraic cycles on quadric sections of cubics in $\PR ^4$]
{\bf Algebraic cycles on quadric sections of cubics in $\PR ^4$ under the action of symplectomorphisms}
\author{V. Guletski\u \i , A. Tikhomirov}

\date{24 March 2014}



\begin{abstract}
\noindent Let $\tau $ be the involution changing the sign of two coordinates in $\PR ^4$. We prove that $\tau $ induces the identity action on the second Chow group of the intersection of a $\tau $-invariant cubic with a $\tau $-invariant quadric hypersurfaces in $\PR ^4$. Let $l_{\tau }$ and $\Pi _{\tau }$ be the $1$- and $2$-dimensional components of the fixed locus of the involution $\tau $. We describe the generalized Prymian associated to the projection of a $\tau $-invariant cubic $\bcC \subset \PR ^4$ from $l_{\tau }$ onto $\Pi _{\tau }$ in terms of the Prymians $\bcP _2$ and $\bcP _3$ associated to the double covers of two irreducible components, of degree $2$ and $3$ respectively, of the reducible discriminant curve. This gives a precise description of the induced action of the involution $\tau $ on the continuous part in the Chow group $CH^2(\bcC )$. The action on the subgroup corresponding to $\bcP _3$ is the identity, and the action on the subgroup corresponding to $\bcP _2$ is the multiplication by $-1$.
\end{abstract}








\maketitle

\section{Introduction}
\label{s-intro}

The aim of this paper is to have yet another look at the correlation between $0$-cycles on algebraic surfaces, on one side, and codimension $2$ algebraic cycles of their $3$-dimensional spreads, on the other. For that purpose we have chosen a fairly concrete model, the intersections of quadric and cubic hypersurfaces invariant under the involution $\tau $ changing the sign of two homogeneous coordinates in $\PR ^4$. Such intersections are $K3$-surfaces, and $\tau $ induces symplectomorphic actions on their second cohomology groups. We will study the induced actions on the second Chow groups of $\tau $-invariant cubics and their $\tau $-invariant quadric sections in $\PR ^4$.

The primary motivation comes from the Bloch-Beilinson conjecture on mixed motives, \cite{Jannsen}, which implies that the induced action of a symplectomorphism of a $K3$-surface on its Chow group of $0$-cycles must be the identity, \cite{Huybrechts}. The first results along this line had been obtained in \cite{Sur les 0-cycles}, where such identity action was proved for quartics in $\PR ^3$ and intersections of $3$ quadrics in $\PR ^5$, covering the cases $d=2$ and $d=4$ in terms of \cite{GeemenSarti}. In \cite{Pedrini} the identity action on the second Chow group was proved for the case $d=1$, and for $K3$-surfaces admitting elliptic pencils with sections.

In the first part of the paper we slightly generalize the method developed in \cite{Sur les 0-cycles} and apply it to prove the identity action of the involution $\tau $ on the second Chow group for intersections of invariant cubics and quadrics in $\PR ^4$ (Theorem \ref{identity}). It should be noted that when the second version of our manuscript was published on the web, we learnt about the paper \cite{HuybrechtsKemeny}, where the identity action of symplectic involutions on the second Chow group was proved in one third of the moduli, and soon after C. Voisin proved the same for all symplectic involutions on $K3$-surfaces, \cite{VoisinNewPaper}. In \cite{Huybrechts2} Huybrechts showed that these results suffice to prove the identity action of symplectomorphisms of any finite order.

In the second part we deal with the induced action of the involution $\tau $ on the second Chow group of a $\tau $-invariant cubic hypersurface $\bcC $ in $\PR ^4$. Our approach is based on the following geometric idea. The set of fixed points of the involution $\tau $ is a union of a line $l_{\tau }$ and a plane $\Pi _{\tau }$ in $\PR ^4$. Projecting $\bcC $ from $l_{\tau }$ onto the plane $\Pi _{\tau }$ we observe that the corresponding discriminant curve splits into a conic $C_2$ and a cubic $C_3$ in $\Pi _{\tau }$. Following \cite{Shokurov}, we construct an isogeny from the generalized Prymian $\bcP $ associated to the double cover of the whole discriminant curve $C_2\cup C_3$ onto the direct product of two Prymians $\bcP _2$ and $\bcP _3$, corresponding to the double covers of $C_2$ and $C_3$ respectively. An enjoyable thing here is that the involution $\tau $ induces the identity action on $\bcP _3$, and it is the multiplication by $-1$ on $\bcP _2$. This gives a complete description of the induced action of the involution $\tau $ on $\bcP $, on the continuous part $A^2(\bcC )$ in the second Chow group of the cubic $\bcC $, as well as on the Hodge pieces of its third cohomology group $H^3(\bcC ,\CC )$. In particular, the action on the subgroup in $A^2(\bcC )$ corresponding to $\bcP _3$ is the identity, and the action on the subgroup corresponding to $\bcP _2$ is the multiplication by $-1$ (Theorem \ref{cubic action}). Pulling back algebraic cycles to the generic fibre of a pencil of quadrics sections we will see that those algebraic cycles which correspond to $\bcP _2$ vanish in the Chow group of the generic fibre.

In a sense, this gives a geometrical reflection, in terms of split discriminant curves and corresponding Prymians, of the behaviour of codimension $2$ algebraic cycles predicted by the Bloch-Beilinson conjectures.

\medskip

{\sc Acknowledgements.} The authors are grateful to Sergey Gorchinskiy for useful suggestions and Claire Voisin for explanation of certain phenomena in algebraic cycles and encouraging interest to this work. 
The research was partially supported by the EPSRC grant EP/I034017/1. The second named author has been financially supported by the Ministry of
Education and Science of Russian Federation, and he also acknowledges the hospitality of the Max Planck Institute for Mathematics in Bonn during the winter 2014.

\section{Notation and terminology}
\label{terminology}

Throughout the paper we work over $\CC $. By default, $H^*(-,A)$ are the Betti cohomology with coefficients $A=\ZZ $, $\QQ $ or $\CC $. The Chow groups will be with coefficients in $\ZZ $. For any quasi-projective variety $X$ over $\CC $, and for any positive integer $q$, let $CH_q(X)$ be the Chow group of dimension $q$ algebraic cycles modulo rational equivalence on $X$. If $X$ is equidimensional of dimension $d$ then the group $CH_q(X)$ will be often denoted by $CH^p(X)$, where $p=d-q$. Our main object of study will be the subgroups $A^p(X)\subset CH^p(X)$ generated by cycles algebraically equivalent to zero on $X$.

Let $X$ be a nonsingular projective complex variety and let $CH^p(X)_{\hom }$ be the kernel of the cycle class homomorphism $cl:CH^p(X)\to H^{2p}(X,\ZZ )$. Each group $H^i(X,\ZZ )$ carries a pure weight $i$ Hodge structure on it. Let $F^p$ be the corresponding decreasing Hodge filtration on the complexified vector space $H^i(X,\CC )$, and let $H^{p,q}(X)$ be the adjoint quotient $(F^{p}/F^{p+1})H^{p+q}(X,\CC )$. The filtration $F^p$ is opposite to the complex conjugate filtration $\bar F^p$, in the sense that $F^p\oplus \overline{F^{q+1}}=H^{p+q}$. Let
  $$
  J^p(X)=H^{2p-1}(X,\CC )/(\im (H^{2p-1}(X,\ZZ ))+F^pH^{2p-1}(X,\CC ))
  $$
be the $p$-th intermediate Jacobian of $X$. Here $\im (H^{2p-1}(X,\ZZ))$ is the image of the natural homomorphism from integral to complex cohomology, i.e. the group $H^{2p-1}(X,\ZZ )$ modulo the torsion subgroup in it. The Poincar\'e duality gives the isomorphism
  $$
  J^p(X)\simeq
  (F^{d-p+1}H^{2d-2p+1}(X,\CC ))^{\vee }/
  \im (H_{2d-2p+1}(X,\ZZ ))\; ,
  $$
where $d$ is the dimension of $X$ and $\im (H_{2d-2p+1}(X,\ZZ ))$ is the image of the integral homology in the dual space of the $d-p+1$-th term of the Hodge filtration via integration of forms over topological chains. Integrating over $2d-2p+1$-dimensional topological chains whose boundaries are homologically trivial $d-p$-cycles gives the Abel-Jacobi homomorphism $AJ:CH^p(X)_{\hom }\to J^p(X)$.

For any nonsingular projective complex variety $Y$ and a cycle $Z$ of codimension $p$ on $Y\times X$ flat over $Y$, we can fix a closed point $y_0$ on $Y$ and define a map from $Y$ to $J^p(X)$ by sending a closed point $y\in Y$ to the image of the algebraically trivial cycle $Z_y-Z_{y_0}$ under the Abel-Jacobi homomorphism $AJ$. This extends to the homomorphism $A_0(Y)\to J^p(X)$, where $A_0(Y)$ is the subgroup generated by algebraically trivial $0$-cycles in $CH_0(Y)$. The image of this homomorphism is a complex subtorus $T_{Y,Z}$ in $J^p(X)$ whose tangent space at $0$ is contained in $H^{p-1,p}(X)$. Let $J_{\alg }^p(X)$ be the maximal subtorus of $J^p(X)$ having this property, so that $J_{\alg }^p(X)$ contains $T_{Y,Z}$ for all possible $Y$ and $Z$. Then $J_{\alg }^p(X)$ is an abelian variety over $\CC $, which is a functor in $X$. The homomorphism $AJ$ sends $A^p(X)$ into $J_{\alg}^p(X)$, so that we obtain the Abel-Jacobi homomorphism $AJ:A^p(X)\to J_{\alg}^p(X)$ on algebraic parts. Notice that the latter homomorphism is expected to be surjective. This is not known in general, but will be the case in the concrete applications below.

All the details on intermediate Jacobians and Abel-Jacobi homomorphisms can be found, for example, in Chapter 12 of \cite{Voisin Book 1}.

Throughout the paper, for an abelian group $A$ and a finite group $G$ acting on $A$ let $N:A\to A$ be the integral averaging operator sending $a\in A$ to $\Sigma _{g\in G}g(a)$, let $A^G$ be the subgroup of $G$-invariant elements in $A$, and let $A^{\sharp }$ be the kernel of the endomorphism $N$. If the group $A$ is divisible and torsion free, then it is a direct sum of $A^G$ and $A^{\sharp }$. Let $\znak $ stay for $G$ or $\sharp $ simultaneously.

Assume now that a finite group $G$ acts by regular automorphisms on $X$. For any $g\in G$ let $g^*$ be the induced automorphism of $H^p(X,A)$, where $A$ is $\ZZ $, $\QQ $ or $\CC $. Each $g^*$ preserves the degrees of differential forms. This is why it preserves $H^{p,q}(X)$ and so the Hodge filtration. Moreover, $g^*$ is compatible with the integration of forms, which gives the automorphism $g^*$ of the complex torus $J^p(X)$ induced by the above automorphisms of complex and integral cohomology groups. By the same reason, $g^*$ is compatible with the corresponding automorphism of $CH^p(X)_{\hom }$ via the Abel-Jacobi homomorphism $AJ$. Since $J_{\alg }^p(X)$ is a functor of $X$, the automorphism $g^*$ gives the automorphism of the abelian variety $J_{\alg }^p(X)$. The norm $N=\sum _{g\in G}g^*$ is an endomorphism of $J_{\alg }^p(X)$ as an abelian variety over $\CC $. Note that $J_{\alg }^p(X)^{\znak }$ is an abelian subvariety in $J^p_{\alg }(X)$ and one has the Abel-Jacobi homomorphisms
  $$
  AJ^{\znak } : A^p(X)^{\znak }\lra J_{\alg}^p(X)^{\znak }
  $$
for $\znak =G$ and $\znak =\sharp $, which will play an important role in what follows.

Next, let $d=2$ and let $\NS (X)$ be the N\'eron-Severy group of the surface $X$. The space $\NS (X)\otimes \QQ $ can be identified with the image $H^2(X,\QQ )_{\alg }$ of the Chow $\QQ $-vector space $CH^1(X)\otimes \QQ $ under the cycle class map to $H^2(X,\QQ )$. The $\cup \; $-product on $H^2(X,\QQ )$ is non-degenerate by the Poincar\'e duality theorem, and it remains non-degenerate after the restriction to $H^2(X,\QQ )_{\alg }$ by the Hodge index theorem. This is why we can consider the orthogonal complement $H^2(X,\QQ )_{\tr }$ to $H^2(X,\QQ )_{\alg }$ with respect to the intersection pairing on $H^2(X,\QQ )$. The group $H^2(X,\QQ )$ is called algebraic if $H^2(X,\QQ )_{\tr }$ is trivial. This is equivalent to say that $p_g=0$, where $p_g=\dim H^{2,0}(X)$ is the geometric genus of the surface $X$. The action of $G$ is compatible with the complex conjugation on the Dolbeault cohomology. This implies that $H^{2,0}(X)^G=0$ if and only if $H^{0,2}(X)^G=0$, and $H^{2,0}(X)^{\sharp }=0$ if and only if $H^{0,2}(X)^{\sharp }=0$. Thus, if $H^{2,0}(X)^G=0$ then $H^2(X,\CC )^G=H^{1,1}(X)^G$, whence the $\QQ $-vector space $H^2(X,\QQ )^G$ is algebraic in the sense that any cohomology class in $H^2(X,\QQ )^G$ comes from the class in $\NS (X)^G\otimes \QQ $ via the cycle class map. Similarly, if $H^{2,0}(X)^{\sharp }=0$ then $H^2(X,\CC )^{\sharp }=H^{1,1}(X)^{\sharp }$, so that the $\QQ $-vector space $H^2(X,\QQ )^{\sharp }$ is algebraic, i.e. a cohomology class in $H^2(X,\QQ )^{\sharp }$ comes from a certain class in the Neron-Severi group $\NS (X)^{\sharp }\otimes \QQ $ via the cycle class map.

In the second half of the paper we will be mainly interested in the case when $G$ is a group of order $2$ generated by an involution $\tau $ acting on a nonsingular projective variety $X$. In particular, if $X$ is a $K3$-surface over $\CC $, $\omega \in H^{2,0}(X)$ a symplectic form on $X$ and $\tau ^*(\omega )=\omega $, then we say that $\tau $ is a symplectomorphism of order $2$ or a Nikulin involution on $X$.

\section{Voisin's theorem}
\label{voisin}

Let $\bcX $ be a nonsingular projective threefold and let $f:\bcX \dasharrow \PR ^1$ be a pencil of surfaces on $\bcX $ with base locus $B$. We will assume that $f$ is nice in the sense that $B$ is nonsingular and that it is the irreducible transversal intersection of two generic members of the pencil. Let $\tilde \bcX \to \bcX $ be the blow up of $\bcX $ at $B$ giving the regular map $\tilde f: \tilde \bcX \to \PR ^1$. The short exact sequence for Chow groups under blowup, \cite[Section 6.7]{Fulton}, yields the isomorphism $A^2(\tilde \bcX )\simeq A^2(\bcX )\oplus A^1(B)$. The Jacobian $\Jac (B)$ of the curve $B$ is a direct summand of the intermediate Jacobian $J^2(\tilde \bcX )$, see \cite{ClemensGriffiths} or \cite{Tjurin}. As taking algebraic parts in intermediate Jacobians is functorial, we obtain the regular surjective morphism of abelian varieties $\epsilon :J_{\alg }^2(\tilde \bcX )\to \Jac (B)$. Assume that a finite group $G$ acts on $\bcX $ fibre-wise. Then $G$ acts also on $B$, and $\epsilon $ induces the regular epimorphism $\epsilon ^{\znak }:J_{\alg }^2(\tilde \bcX )^{\znak }\to \Jac (B)^{\znak }$
which will be used later. For any $t\in \PR ^1$ let $X_t$ be the fibre of the regular morphism $\tilde f$. We impose the following assumption on the pencil $f:\bcX \dasharrow \PR ^1$.

\medskip

\begin{itemize}

\item[(A)]{}
{\it $\exists $ a Zariski open $U\subset \PR ^1$, such that $X_t$ is nonsingular, $H^2(X_t,\QQ )^{\znak }$ is algebraic and $H^1(X_t,\QQ )^{\znak }=0$ for any $t\in U$.}

\end{itemize}

\medskip

For any $t\in \PR ^1$ let $j_t:B\to X_t$ be the closed embedding of the base locus into the fibre. Since $B$ is a Cartier divisor in $X_t$, the embedding $j_t$ induces the Gysin homomorphism $j_t^*:CH^1(X_t)\to CH^1(B)$, see \cite[Section 2.6]{Fulton}.

\begin{theorem}
\label{Voisin's theorem}
Under Assumption (A), the group $A^1(B)^{\znak }$ is contained in the image of the natural homomorphism $\oplus _{t\in \PR ^1}CH^1(X_t)^{\znak }\to CH^1(B)^{\znak }$ induced by the homomorphisms $j_t^*$.
\end{theorem}

\begin{pf}
Let $U$ be as in (A), let $\eta $ be the generic and $\bar \eta $ the geometric generic point of $\PR ^1$. Consider the relative Picard scheme $\bcP \to U$ of the pull-back $\tilde f_U:\tilde \bcX _U\to U$ of $\tilde f$ to $U$. Its fibre $\bcP _{\bar \eta }$ is the Picard scheme of the fibre $X_{\bar \eta }$ and, by Assumption (A), $\bcP _{\bar \eta }^{\znak }$ is isomorphic to the N\'eron-Severi group $\NS (X_{\bar \eta })^{\znak }$. The group $G$ acts in the fibres of the structural morphism from $\bcP $ onto $U$. Since $\NS (X_{\bar \eta })$ is finitely generated, we choose a finite number of points $P_1,\dots ,P_n$ in $\bcP _{\bar \eta }^{\znak }$ generating the group $\NS (X_{\bar \eta })^{\znak }$. For each $P_i$ let $W_i\to U$ be a finite morphism onto the curve $U$, such that the curve $W_i$ is nonsingular and the residue field of the scheme $\bcP _{\eta }$ at $P_i$ is $\CC (W_i)$. Let $W^{\znak }$ be the union of the curves $W_i$ with the structural morphism $g^{\znak }:W^{\znak }\to U$. Then, for any $t\in U$, the fibre $W^{\znak }_t$ generates the $\QQ $-vector space $\NS (X_t)_{\QQ }^{\znak }$, which in turn is isomorphic to $H^2(X_t,\QQ )^{\znak }$ by the assumption imposed on $f$.

As $U$ is locally contractible, the stalk of $R^0g^{\znak }_*\QQ $ at $t\in U$ is $H^0(W^{\znak }_t,\QQ )$ and the stalk of $R^2\tilde f_*\QQ $ at $t$ is $H^2(X_t,\QQ )$. Shrinking $U$ if necessary we can assume that $g^{\znak }$ and $f_U$ are smooth. Trivializing these morphisms in complex topology one can define the sheaf $(R^2(\tilde f_U)_*\, \QQ )^{\znak }$ and the surjective morphism
  $$
  \alpha : R^0g^{\znak }_*\QQ \lra
  (R^2(\tilde f_U)_*\, \QQ )^{\znak }
  $$
of sheaves on $U$, such that for any $t\in U$ and $P\in W_t^{\znak }$ the local homomorphism $\alpha _t$ sends $P_i$ to the corresponding element in $\NS (X_t)$. Since $U$ is homotopy equivalent to a $1$-dimensional $CW$-complex, it has cohomological dimension $1$, so that the induced homomorphism
  $$
  \alpha _* : H^1(U,R^0g^{\znak }_*\QQ )\to
  H^1(U,(R^2(\tilde f_U)_*\, \QQ )^{\znak })
  $$
is surjective.

Let $\bcY _i$ be the fibred product $W_i\times _U\tilde \bcX _U$ and let $\bcS _i\to W_i$ be the Picard scheme of the relative scheme $\bcY _i\to W_i$. Each point $P_i$ is rational over the field $\CC (W_i)$ bringing a section of the morphism $\bcS _i\to W_i$ over some possibly smaller Zariski open subset $V_i$ in $W_i$. This section in turn induces the section of the morphism $\bcY _i\times _{W_i}\bcS \to \bcY _i$ over $(\bcY _i)_{V_i}$. Let $\bcD _i$ be the pull-back of the Poincar\'e divisor on the scheme $\bcY _i\times _{W_i}\bcS _i$ to $(\bcY _i)_{V_i}$, and let $\bcD $ be the union of $\bcD _i$'s. Using appropriate nonsingular compactifications $\bar W_i$ and considering their union $\bar W^{\znak }$ we can also work with the morphism $\bar g^{\znak }:\bar W^{\znak }\to \PR ^1$, such that $g^{\znak }$ is the pull-back of $\bar g^{\znak }$ under the embedding of $U$ into $\PR ^1$. Let $\bar \bcD _i$ be the closure of $\bcD _i$ in the fibred product $\bar W^{\znak }\times \tilde \bcX $ over $\Spec (\CC )$ and $\bar \bcD $ be the union of $\bar \bcD _i$'s. Then $\bar \bcD $ is an algebraic cycle of codimension $2$ in the fourfold $\bar W^{\znak }\times \tilde \bcX $. Let
  $$
  \bar \beta =
  cl(\bar \bcD )\in H^4(\bar W^{\znak }\times \tilde \bcX ,\ZZ )
  $$
be the cohomology class of $\bar \bcD $. Passing to rational coefficients in cohomology groups, the $(1,3)$-K\"unneth component $\bar \beta (1,3)$ of $\bar \beta $ induces the homomorphism of Hodge structures
  $$
  \bar \beta (1,3)_* : H^1(\bar W^{\znak },\QQ )\to
  H^3(\tilde \bcX ,\QQ )\; .
  $$

Let $\bar \beta (1,3)_{i,*}$ be its restriction on $H^1(\bar W^{\znak }_i,\ZZ )$ and let $J_i$ be the Jacobian of the curve $\bar W_i$. Then $\bar \beta (1,3)_{i,*}$ induces the regular morphism
  $$
  \bar \beta (1,3)_{i,*} : J_i\lra J^2(\tilde \bcX )\; ,
  $$
which factorizes through $J_{\alg }^2(\tilde \bcX )$. For any two points $P$ and $P_0$ on $W_i$, let $D_P$ and $D_{P_0}$ be the divisors on $X_t$ and $X_{t_0}$ respectively, whose cohomology classes correspond to $P$ and $P_0$ as points in the fibres of the morphism $\bcP \to U$. Then $\bar \beta (1,3)_{i,*}([P-P_0])=AJ[D_P-D_{P_0}]$ in $J_{\alg }^2(\tilde \bcX )$, see Theorems 12.4 and 12.17 in \cite{Voisin Book 1}.

Next, $R^1g^{\znak }_*\QQ =0$ as the fibres $W_t^{\znak }$ are $0$-dimensional, and $(R^3(\tilde f_U)_*\QQ )^{\znak }=0$ by Assumption (A). Since the corresponding Leray spectral sequences $E_2$-degenerate, we obtain the isomorphisms
  $$
  H^1(U,R^0g^{\znak }_*\QQ )\simeq H^1(W^{\znak },\QQ )
  \quad \hbox{and}\quad
  H^1(U,(R^2(\tilde f_U)_*\QQ )^{\znak })\simeq H^3(\tilde \bcX _U,\QQ )^{\znak }\; .
  $$
Due to the above description of the action of $\bar \beta (1,3)_{i,*}$ on $[P-P_0]$, one can check that the diagram
    $$
    \xymatrix{
    H^1(\bar W^{\znak },\QQ )
    \ar[dd]^-{r_1} \ar[rr]^-{\bar \beta (1,3)_*} & &
    H^3(\tilde \bcX ,\QQ )^{\znak } \ar[dd]^-{r_2} \\ \\
    H^1(W^{\znak },\QQ )
    \ar[rr]^-{\zeta _*} & &
    H^3(\tilde \bcX _U,\QQ )^{\znak }
    }
    $$

\medskip

\noindent commutes, where $\zeta _*$ is the modification of $\alpha _*$ by means of the above two isomorphisms coming from the Leray spectral sequences, and $r_1$, $r_2$ are the restriction homomorphisms on cohomology groups.

By Deligne's results, the cohomology groups at the bottom possess mixed Hodge structures with weights
  $$
  W_0H^1(W^{\znak },\QQ )=\im (r_1)\qquad \hbox{and}\qquad W_0H^3(\tilde f^{-1}(U),\QQ )^{\znak }=\im (r_2)\; ,
  $$
see \cite[Section 11.1.4]{Voisin Book 2}. Morphisms between mixed Hodge structures are strict with respect to both Hodge and weight filtrations, see \cite{HodgeTheoryII}. Since $\zeta _*$ is a morphism of mixed Hodge structures and it is surjective, we obtain that
  $$
  \zeta _*(W_0H^1(W^{\znak },\QQ ))=
  W_0H^3(\tilde \bcX _U,\QQ )^{\znak }\; .
  $$
This gives that $\im (r_2)=\im (r_2\circ \bar \beta (1,3)_*)$. Then $H^3(\tilde \bcX ,\ZZ )^{\znak }$ is generated by $\ker (r_2)$ and $\im (\beta (1,3)_*)$.

On the other hand, $\ker (r_2)$ is generated by the images of the homomorphisms
  $$
  (i'_t)_*:H^1(X'_t,\QQ )^{\znak }\to
  H^3(\tilde \bcX ,\QQ )^{\znak }\; ,
  $$
where $t\in \PR ^1\smallsetminus U$, $X'_t$ is the resolution of singularities of $X_t$ and $i'_t$ is the composition of the desingularization $X_t'\to X_t$ with the closed embedding $i_t:X_t\to \bcX $.

All these things together give that the homomorphism
  $$
  \theta : H^1(\bar W^{\znak },\QQ )\oplus (\oplus _{t\in \PR ^1\smallsetminus U}H^1(X'_t,\QQ )^{\znak })\to H^3(\tilde \bcX ,\QQ )^{\znak }\; ,
  $$
induced by the homomorphisms $\bar \beta (1,3)_*$ and $(i_t)_*$, is surjective.

Since $\theta $ is a homomorphism of polarized Hodge structures, it induces a surjective homomorphism of the corresponding abelian varieties
  $$
  \rho ^{\znak }:J^{\znak }\oplus (\oplus _{t\in \PR ^1\smallsetminus U}
  (\bcP '_{t,\, 0})^{\znak })\to J_{\alg }^2(\tilde \bcX )^{\znak }\; ,
  $$
where $J^{\znak }$ is the union of $J_i$'s and $\bcP '_{t,0}$ is the component of $0$ in the Picard scheme $\bcP '_t$ of the surface $X'_t$.

Now, let $\tilde B$ be the exceptional divisor of the blow-up $\tilde \bcX \to \bcX $, $p:\tilde B\to B$ the projection and $e:\tilde B\to \tilde \bcX $ the embedding of $\tilde B$ into $\tilde \bcX $. Let $\varepsilon ^{\znak }:CH^2(\tilde \bcX )^{\znak }\to CH^1(B)^{\znak }$ be the composition of the pull-back $(e^*)^{\znak }:CH^2(\tilde \bcX )^{\znak }\to CH^2(\tilde B)^{\znak }$ and push-forward $p_*^{\, \znak }:CH^2(\tilde B)^{\znak }\to CH^1(B)^{\znak }$, induced by $e$ and $p$ respectively. For each $t\in \PR ^1$ let $\tilde i_t$ be the closed embedding of the fibre $X_t$ into $\tilde \bcX $, and let $(\tilde i_t)_*^{\, \znak }:CH^1(X_t)^{\znak }\to CH^2(\tilde \bcX )^{\znak }$ be the push-forward homomorphism induced by the closed embedding $\tilde i_t$.

The above surjective homomorphisms $\rho ^{\znak }$ and $\epsilon ^{\znak }$, constructed on the level of Jacobians, guarantee that, on the level of Chow groups, $A^1(B)^{\znak }$ is contained in the image of the sum
  $$
  \oplus _{t\in \PR ^1}
  (\varepsilon ^{\znak }\circ (\tilde i_t)_*^{\, \znak }):
  \oplus _{t\in \PR ^1}CH^1(X_t)^{\znak }\to
  CH^1(B)^{\znak }
  $$
A straightforward verification shows that, for each $t\in \PR ^1$, the composition $\varepsilon ^{\znak }\circ (\tilde i_t)_*^{\, \znak }$ is the restriction of the pull-back homomorphism $j_t^*:CH^1(X_t)\to CH^1(B)$ onto the $\znak $-parts of the Chow groups. This completes the proof of Theorem \ref{Voisin's theorem}.
\end{pf}

\begin{remark}
\label{move}
{\rm Theorem \ref{Voisin's theorem} can be strengthened by saying that $A^1(B)^{\znak }$ is contained in the image of the homomorphism $\oplus _{t\in \PR ^1\smallsetminus Z}CH^1(X_t)^{\znak }\to CH^1(B)^{\znak }$, where $Z$ is a finite subset in $U$. This is because we can always move a zero-cycle in its class modulo rational equivalence on $\bar W^{\znak }$.
}
\end{remark}

\section{The $\tau $-action on $CH^2(S)$}
\label{idaction}

Symplectomorphisms over $\CC $ have order $\leq 8$, see \cite{Nikulin}. If $\tau $ is a Nikulin involution, then $\rho \geq 9$, where $\rho $ is the rank of the N\'eron-Severi group $NS(X)$, see \cite{GeemenSarti}. Assume that $\rho =9$ and let $L$ be a generator of the orthogonal complement of the lattice $E_8(-2)$ in $NS(X)$ whose self-intersection is $2d$, for some positive integer $d$. Let $\Gamma $ be a direct sum of $\ZZ L$ and $E_8(-2)$, if the integer $d$ is odd, or the unique even lattice containing $\ZZ L\oplus E_8(-2)$ as a sublattice of index $2$, if $d$ is even. For each $\Gamma $ there exists a $K3$-surface $X$ with a Nikulin involution and $\rho =9$, such that $NS(X)\simeq \Gamma $, and all such surfaces are parametrized by a coarse moduli space of dimension $11$, see \cite{GeemenSarti}, Proposition 2.3.

Let $S_0$ be a $K3$-surface over $\CC $ with a Nikulin involution $\tau $, such that $\rho =9$ and $d=3$. In this case the generator $L$ gives the regular embedding $\phi _L:S_0\to \PR ^4$, which identifies $S_0$ with the complete intersection of nonsingular cubic $\bcC _0$ and quadric $\bcQ _0$ in $\PR ^4$. The involution $\tau $ extends to the involution $\tau _{\PR ^4}$ on the whole projective space $\PR ^4$. In suitable coordinates, $\tau _{\PR ^4}$ sends $(x_0:x_1:x_2:x_3:x_4)$ to $(-x_0:-x_1:x_2:x_3:x_4)$. The cubic $\bcC _0$ and quadric $\bcQ _0$ are both invariant under $\tau _{\PR ^4}$. Vice versa, if $\bcC _0$ and $\bcQ _0$ are general nonsingular cubic and quadric in $\PR ^4$, both invariant under the involution $\tau _{\PR ^4}$, and such that their intersection $S_0=\bcC _0\cap \bcQ _0$ is nonsingular, then $S_0$ is a $K3$-surface with the Nikulin involution $\tau =\tau _{\PR ^4}|_{S_0}$, see Section 3.3 in \cite{GeemenSarti}.

For short, we will write $\tau $ for $\tau _{\PR ^4}$ and for the involution on $S_0$ simultaneously. The fixed locus of $\tau $ on $\PR ^4$ is the disjoint union of the line $l_{\tau }$ and the plane $\Pi _{\tau }$ given by the equations
  \begin{equation}
  \label{tau-invariant}
  l_{\tau } : x_2 = x_3 = x_4 = 0\; ,\qquad
  \hbox{and}\qquad \Pi _{\tau } : x_0 = x_1 = 0\; .
  \end{equation}
Let $V$ be a vector space, such that $\PR ^4=\PR (V)$. In coordinate-free terms, $\tau $ lifts to the involution $\tau :V\to V$ which induces two involutions $\tau _i:\Sym ^iV^{\vee }\to \Sym ^iV^{\vee }$, where $i=2,3$, $\Sym ^i$ is the $i$-th symmetric power of a vector space over $k$, and $V^{\vee }$ is the $k$-vector space dual to $V$. Consider the subspaces
  $$
  (\Sym ^iV^{\vee })_+ = \{ F\in \Sym ^iV^{\vee }\; |\; \tau _i(F)=F\} \; .
  $$

\medskip

\noindent for $i=2,3$. Any $F\in (\Sym ^2V^{\vee })_+$ has the shape
  \begin{equation}
  \label{S2+}
  \alpha_{00}x_0^2+\alpha_{11}x_1^2+\alpha_{01}x_0x_1+
  f_2(x_2,x_3,x_4)
  \end{equation}
and $\Phi \in (\Sym ^3V^{\vee })_+$ has the shape
  \begin{equation}
  \label{S3+}
  l_{00}(x_2,x_3,x_4)x_0^2+l_{11}(x_2,x_3,x_4)x_1^2+l_{01}(x_2,x_3,x_4)x_0x_1 +f_3(x_2,x_3,x_4)\; ,
  \end{equation}
where $\alpha _{ij}$ are constants, $l_{ij}$ are linear forms, $f_2$ and $f_3$ are homogeneous polynomials of degree $2$ and $3$ respectively. If
  $$
  \bcL _2 = \PR ((\Sym ^2V^{\vee })_+)\quad \hbox{and}\quad
  \bcL _3 = \PR ((\Sym ^3V^{\vee })_+)\; ,
  $$
then $\bcQ _0\in \bcL _2\subset |\bcO (2)|$ and $\bcC_0\in \bcL _3\subset |\bcO (3)|$.

The explicit formulae (\ref{tau-invariant}) and (\ref{S3+}) show that any cubic $\bcC \in \bcL _3$ contains the line $l_{\tau }$. From (\ref{S3+}) it follows that $\bcL _3$ is spanned by the subsystem
  $$
  \bcL _{3,i}=\PR (V_i)\; ,\qquad i=1,2\; ,
  $$
where $V_1$ is the subspace of forms in $\Sym ^3V^{\vee }$ of the shape $$
  l_{00}(x_2,x_3,x_4)x_0^2+l_{11}(x_2, x_3, x_4)x_1^2+
  l_{01}(x_2, x_3, x_4)x_0x_1
  $$
and $V_2$ is the subspace of forms in $\Sym ^3V^{\vee }$ of the shape
  $$
  f_3(x_2,x_3,x_4)\; ,
  $$
and the forms $l_{ij}$, $f_3$ are those described above. The subgroup
  $$
  G=\{ g\in \PGL (5)\; \, |\; \, g(l_{\tau })=
  l_{\tau }\; , g(\Pi _{\tau })=\Pi _{\tau }\}
  $$
acts transitively on the set $\PR ^4\smallsetminus (l_{\tau }\sqcup \Pi _{\tau })$. It also acts naturally on the linear system $|\bcO (3)|$ and fixes the subspaces $\bcL _{3,1}$ and $\bcL _{3,2}$ in it. Then $\bcL _3$ is fixed under the $G$-action too. Let $P_0=(1:1:1:1:1)$ be the point in $\PR ^4$. From (\ref{S3+}) we have that the subsystem $\{ \bcC \in \bcL _3\; |\; P_0\in \bcC \} $ is a proper hyperplane in $\bcL _3$. Since $G$ acts transitively on $\PR ^4\smallsetminus(l_{\tau }\sqcup \Pi _{\tau })$, it follows that, for any point $P\in \PR ^4\smallsetminus (l_{\tau }\sqcup \Pi _{\tau })$ the subsystem $\{ \bcC \in \bcL _3\; |\; P\in \bcC \} $ is also a proper hyperplane of $\bcL $. This shows that the line $l_{\tau }$ is the base locus of the linear system $\bcL _3$.

Consider the linear system
  $$
  \bcM _3=\{ \bcC \in \bcL_3 \; |\; S_0 \subset \bcC\}
  $$
and its subset $\bcN _3$ of the unions $\bcQ _0 \cup H\in \bcM _3$, such that $H$ is a hyperplane containing the line $l_\tau $. Then
  $$
  \bcN _3\simeq \PR ^2\qquad \hbox{and}\qquad
  \bcM _3= \spn _{\bcL _3}(\bcN _3,\bcC _0)
  \simeq \PR ^3
  $$
We also need the linear system
  $$
  \Sigma = \{ B_1\subset S_0 \; | \; B_1=S_0 \cap \bcC _1\; , \; \; \bcC _1\in \bcL _3\smallsetminus \bcM _3\}
  $$
of curves on $S_0$ cut out by cubics from $\bcL _3\smallsetminus \bcM _3$. As $\bcL _3\cong \PR ^{18}$, it follows that $\Sigma \simeq \PR ^{14}$. For any point $P$ in $S_0$, the set
  $$
  \Sigma _P=\{B_1\in\Sigma \; | \; P\in B_1\}
  $$
is a hyperplane in $\Sigma $ if $P\not \in S_0\cap l_{\tau }$, and $\Sigma _P=\Sigma $ otherwise. In both cases, as $\dim (\Sigma )=14$, for any two distinct points $P$ and $Q$ in $S_0$ we have that
  $$
  \dim (\Sigma _P\cap \Sigma _Q)\geq 12\; .
  $$

\begin{lemma}
\label{ermak}
For a general choice of $\bcC _0$ and $\bcQ _0$, there is a nonempty Zariski open subset $N_0$ in $S_0$, such that, if $P_0$ is a point in $N_0$, one can find a nonempty Zariski open subset $U_0$ in $S_0$ having the property that for each point $P$ in $U_0$ there exists a nonsingular curve $B_1\in \Sigma $ passing through $P$ and $P_0$ on $S_0$.
\end{lemma}

\begin{pf}
For a general $\bcQ _0$ in $\bcL _2$, the quadric $\bcQ _0$ intersects $l_{\tau }$ at two distinct points, say $P_+$ and $P_-$. As $l_{\tau }$ is the base locus of $\bcL _3$, the union $S_0\cap l_{\tau }=P_+\sqcup P_-$ is the base locus of the linear system $\Sigma $ on $S_0$. Let $V_+$ be the set of all triples $(\bcC _0,\bcC _1,\bcQ )\in \bcL _3\times \bcL _3\times \bcL _2$, such that $\bcC _0\neq \bcC _1$, the set $B_1=\bcC _0\cap \bcC _1\cap \bcQ _0$ has dimension $1$ in a Zariski open neighbourhood of the point $P_+$, and $B_1$ is nonsingular at $P_+$. Then $V_+$ is a Zariski open subset in $\bcL _3\times \bcL _3\times \bcL _2$. In appropriate coordinates, $P_+=(1:0:0:0:0)$ and $P_-=(0:1:0:0:0)$. Then the triple $\bcC _0=\{ x_2x_0^2=0\} $, $\bcC _1=\{ x_3x_0^2=0\} $, $\bcQ _0=\{ x_0x_1=0\} $ is in $V_+$, whence $V_+\neq \emptyset $. Similarly, one can construct the nonempty open subset $V_-$ in $\bcL _3\times \bcL _3\times \bcL _2$, regarding the point $P_-$. Joint with Bertini's theorem, this gives that, for a general choice of $\bcC _0$ and $\bcQ _0$, there is a nonempty Zariski open subset $V$ in $\Sigma $, such that each curve $B_1\in V$ is nonsingular.

The set $T=\{ (P,Q)\in S_0\times S_0\, |\, \dim (\Sigma _P\cap \Sigma _Q)=12\} $ is Zariski open, and hence irreducible, in $S_0\times S_0$. The set $Z=\{ (P,Q,B_1)\in T\times \Sigma \, |\, P,Q\in B_1\} $ is Zariski closed in $T\times \Sigma $. The projection $\pi :Z\to T$ is surjective. Since $\pi $ is a $\PR ^{12}$-bundle over $T$ and $T$ is irreducible, $Z$ is irreducible too. Let $B_1$ be any curve in $V$ and $P$ be any point on $B_1$, not equal to $P_0$ or $P_1$. The set $F_P=\{ Q\in S_0\, |\, \dim (\Sigma _P\cap \Sigma _Q)>12\} $ is at most finite. If $Q\in B_1\smallsetminus F_P$ then $(P,Q,B_1)\in Z$ and $(P,Q,B_1)\in S_0\times S_0\times V$. Therefore, $W=Z\cap (S_0\times S_0\times V)$ is a nonempty Zariski open subset in the irreducible quasi-projective variety $Z$.

As $\pi $ is surjective, the Zariski closure of the set $\pi (W)$ is $T$. Since, moreover, the image $\pi (W)$ is constructible, it contains a subset $T_0$, which is open and dense in $T$. Then $T_0$ is Zariski open and dense also in $S_0\times S_0$. It follows that the image of $T_0$ under the projection of $S_0\times S_0$ onto the second factor contains a nonempty Zariski open subset $N_0$. For any point $P_0\in N_0$ let then $U_0$ be the image of the set $T_0\cap (S_0\times \{ P_0\} )$ under the projection of $S_0\times S_0$ onto the first factor.
\end{pf}

\begin{theorem}
\label{identity}
Let $S$ be a general nonsingular complete intersection of cubic and quadric hypersurfaces, both invariant under the involution $\tau $ in $\PR ^4$. Then the action $\tau ^*:CH^2(S)\to CH^2(S)$ is the identity.
\end{theorem}

\begin{pf}
In the above terms, $S=S_0$ is the intersection of $\tau $-invariant nonsingular cubic $\bcC _0$ and quadric $\bcQ _0$ in $\PR ^4$. Let $N_0$ be the nonempty Zariski open subset in $S_0$, coming from Lemma \ref{ermak}, and let $P_0$ be a point in $N_0$. As the action of $\tau ^*$ does not change the degree of $0$-cycles on $S$, to prove the theorem all we need to show is that, for any point $P$ on $S_0$, the cycle class $[P-P_0]$ is $\tau $-invariant. Let $U_0$ be the nonempty Zariski open subset in $S_0$, depending on $P_0$, also as in Lemma \ref{ermak}. By the Chow moving lemma, one can assume that $P\in U_0$. Then, by Lemma \ref{ermak}, there exists a cubic $\bcC _1$ in $\bcL _3\smallsetminus \bcM _3$, such that $\bcC _1$ passes through $P$ and $P_0$, and the curve $B_1=S_0\cap \bcC _1$ is nonsingular. Let $\bcQ =\bcQ _0$, and let $f:\bcQ \dasharrow \PR ^1$ be the pencil of $K3$-surfaces obtained by restricting the pencil $|\bcC _t|_{t\in\PR ^1}$, spanned by $\bcC _0$ and $\bcC _1$, onto the quadric $\bcQ $. The nonsingular curve $B=B_1$ is the base locus of the pencil $f$ and $P,P_0\in B$.

Now, for any $t\in \PR ^1$ let $S_t=\bcC _t\cap \bcQ $, let $i_t:S_t\to \bcQ $ be the corresponding closed embedding, and let $j_t:B\to S_t$ be the closed embedding of the base locus into the fibre (without loss of generality, we may think of $S_0$ as the fibre over the point $t=0$). Let $\alpha $ be the class of $P-P_0$ in $A^1(B)$ and let $\alpha ^{\sharp }=\alpha -\tau ^*(\alpha )\in A^1(B)^{\sharp }$. To prove the theorem we need to show that ${j_0}_*(\alpha ^{\sharp })$ vanishes.

The Assumption (A) is satisfied for the pencil $f$ and $\znak =\sharp $. By Theorem \ref{Voisin's theorem} and Remark \ref{move}, the cycle class $\alpha ^{\sharp }$ is a sum of cycle classes of type $j_t^*(\alpha ^{\sharp }_t)$, where $\alpha ^{\sharp }_t\in CH^1(S_t)^{\sharp }$ and $t\neq 0$. For each $t\in \PR ^1$, such that $t\neq 0$, we have the Cartesian square
    $$
    \xymatrix{
    B \ar[dd]_-{j_t} \ar[rr]^-{j_0} & & S_0 \ar[dd]^-{i_0} \\ \\
    S_t \ar[rr]^-{i_t} & & \bcQ
    }
    $$
It consists of four closed embeddings, each of which is an embedding of a Cartier divisor into the target variety. This is why all the embeddings are regular, whence ${j_0}_*\circ j_t^*=i_0^*\circ {i_t}_*$, as homomorphisms from $CH^1(S_t)$ to $CH^2(S_0)$, see \cite[Section 6.2]{Fulton}. Since $\alpha ^{\sharp }$ is a sum of cycle classes $j_t^*(\alpha ^{\sharp }_t)$, it follows that ${j_0}_*(\alpha ^{\sharp })=i_0^*(\delta ^{\sharp })$ for some $\delta ^{\sharp }$ in $CH^2(\bcQ )^{\sharp }$. Since $\bcQ $ is a $3$-dimensional quadric hypersurface in $\PR ^4$, the group $CH^2(\bcQ )$ is isomorphic to $\ZZ $, with the generator represented by class of a line $L$ in $\bcQ $. As the line $\tau (L)$ is rationally equivalent to $L$ on $\bcQ $, the group $CH^2(\bcQ )^{\sharp }$ vanishes. Therefore, $\delta ^{\sharp }=0$ and hence ${j_0}_*(\alpha ^{\sharp })=0$. This finishes the proof of Theorem \ref{identity}.
\end{pf}

\section{The $\tau $-action on $A^2(\bcC )$}

Let $\bcC $ be a general nonsingular cubic from $\bcL _3$ and consider the linear projection of $\PR ^4$ onto $\Pi _{\tau }$ from the line $l_{\tau }\subset \bcC $. Restricting the projection onto $\bcC $ we get a rational map $p:\bcC \dashrightarrow \Pi _{\tau }$. Blowing up $\bcC $ at the indeterminacy locus $l_{\tau }$ we obtain the conic bundle
  $$
  \hat p:\hat \bcC \to \Pi _{\tau }\; .
  $$
Let
  $$
  C\subset \Pi _{\tau }
  $$
be the discriminant curve of $\hat p$. This is an algebraic curve of degree $5$ in $\Pi _{\tau }$. Let also $\bcF $ be the Fano surface of lines on $\bcC $. Following \cite{MurreCubics}, we look at the Zariski closed subset $\bcF _0$ (respectively, $\bcF '_0$) in $\bcF $ generated by lines $l$, such that for $l$ there exists a plane $\Pi $ in $\PR ^4$ with $\bcC \cdot \Pi = 2l + l'$ (respectively, $\bcC \cdot \Pi = l + 2l'$). In loc. cit. the cubic is projected from a line belonging neither to $\bcF _0$ nor to $\bcF _0'$, so that the discriminant curve is irreducible. In our case the line $l_{\tau }$ is not in $\bcF _0$, but still is an element of $\bcF _0'$. This has the effect that the discriminant curve $C$ is reducible and consists of $2$ irreducible components,
  $$
  C=C_2\cup C_3\; .
  $$
Here $C_2$ is the conic defined by the equation
  $$
  4l_{00}(x_2,x_3,x_4)l_{11}(x_2,x_3,x_4)-
  l_{01}(x_2,x_3,x_4)^2=0\; ,
  $$
and $C_3$ is the cubic defined by the equation
  $$
  f_3(x_2,x_3,x_4)=0
  $$
in $\Pi _{\tau }$, where $l_{00}$, $l_{11}$, $l_{01}$ and $f_3$ are as in (\ref{S3+}). For the general choice of $l_{00}$, $l_{11}$, $l_{01}$ and $f_3$ the curves $C_2$ and $C_3$ are nonsingular and intersect each other transversally at $6$ distinct points in $\Pi _{\tau }$.

For any point $P$ on $\Pi _{\tau }$ the span $\Pi _P$ of $P$ and $l_{\tau }$ intersects the cubic $\bcC $ along the line $l_{\tau }$ and a conic $C_P$. Since $P$, as well as the points of $l_{\tau }$ are fixed under the involution $\tau $, the involution acts in the fibres of the morphism $\hat p$. If $P\in C$ then the conic $C_P$ splits into two lines
  $$
  l_P^+\qquad \hbox{and}\qquad l_P^-\; ,
  $$
which coincide when $P$ belongs to $C_2\cap C_3$.

Let $C'_i$ be the curve generated by the lines $l$ in the Grassmannian $\Gr (2,5)$, such that $p(l\smallsetminus l_{\tau })\in C_i$, for $i=2,3$, and let
  $$
  C'=C'_2\cup C'_3\; .
  $$
Then $C'$ is a double cover of the curve $C$ which induces the involution $\iota :C'\to C'$ by transposing the lines $l_P^+$ and $l_P^-$ sitting over the points $P$ in $C$. The double covers $C_2'\to C_2$ and $C_3'\to C_3$ are ramified over the six points of intersection of $C_2$ and $C_3$, having no ramification in other points of $C'$. Each point $P\in C_2\cap C_3$ is a double point of $C$, and if $P'\in C'_2\cap C'_3$ sits over $P$, then $P'$ is a double point of $C'$, see 1.5.3 in \cite{Beauville}. This is why the curves $C'_2$ and $C'_3$ are nonsingular. The Hurwitz formula shows that the genera of the curves $C_2'$ and $C_3'$ are $2$ and $4$ respectively.

The above involution $\iota $ acts component-wise, which gives two involutions $\iota _2$ on $C_2'$ and $\iota _3$ on $C_3'$. The curve $C'$ has only double point singularities lying over the six double points in $C_2\cap C_3$. These six points on $C'$ are the fixed points of the involution $\iota $. Then $(C',\iota )$ is a Beauville pair, i.e. it satisfies the condition (B) on page 100 in \cite{Shokurov}. Let
  $$
  \bcP =\ker (\Nm )^0=(\id -\iota ^*)\Pic (C')
  $$
be the generalized Prymian in the sense of Beauville, see \cite{Beauville} or \cite{Shokurov}. Notice that $\bcP $ is a principally polarized abelian variety over $\CC $, loc.cit. Any closed point $P$ on $C'_2$ or $C'_3$ gives the line $L_P$ on $\bcC $. By Beauville's result, \cite[Section 3.6]{Beauville}, the correspondence $P\mapsto L_P$ induces the isomorphism
  $$
  \bcP \simeq A^2(\bcC )\; .
  $$

Let also $J_i$ be the Jacobian of the curve $C'_i$, for $i\in \{ 2,3\} $. The involutions $\iota _i$ give the Primians
  $$
  \bcP _i = (\id -\iota _i^*)J_i\; .
  $$
Since the genus of $C_2$ is zero, $\bcP _2$ coincides with $J_2$, which is an abelian surface over $\CC $.

Let $N_6\to \dots \to N_1\to C$ be the chain of six subsequent normalizations of the curve $C$ at the six double points of $C$, and let $N'_6\to \dots \to N'_1\to C'$ be the chain of the corresponding six normalizations of the curve $C'$ at the double points lying over the double points of $C$. Each curve $N'_{i+1}$ inherits an involution from $N'_i$, and the corresponding generalized Prymians $\bcR _i$, $i=1,\dots ,6$, are abelian varieties by Theorem 3.5 in \cite{Shokurov}. Moreover, the involutions on the curves $N'_i$ satisfy the the condition (F) of Proposition 3.9 in \cite{Shokurov}, and so give the tower of isogenies $\bcP \to \bcR _1\to \dots \to \bcR _6$ by Lemma 3.15 in loc.cit. The curve $N'_6$ is the disjoint union of the curves $C'_2$ and $C'_3$, and the restrictions of the induced involution on $N'_6$ on the connected components $C'_2$ and $C'_3$ coincide with the involutions $\iota _2$ and $\iota _3$ respectively. Therefore,
  $$
  \bcR _6=\bcP _2\oplus \bcP _3\; ,
  $$
and we obtain the isogeny
  $$
  \Lambda :\bcP \to \bcP _2\oplus \bcP _3\; .
  $$

Notice that by the same Lemma 3.15 in Shokurov's paper, the first isogeny $\bcP \to \bcR _1$ is an isomorphism because the Beauville pair $(C',\iota )$ satisfies the condition (B) and the last isogeny $\bcR _5\to \bcR _6$ is an isomorphism because $N'_6$ is disconnected. Each of the rest $4$ isogenies is of degree $2$, loc.cit. Then the total isogeny $\Lambda $ is of degree $2^4$.

Since $\bcC $ is unirational, there exists the classically known rational dominant morphism $\PR ^3\dasharrow \bcC $, see \cite{ClemensGriffiths}, Appendix B. Resolving its indeterminacy, we get the dominant regular morphism $\hat \PR ^3\to \bcC $.  Then $\hat \PR ^3$ is balanced by Prop. 1.2 in \cite{BV}, and $\bcC $ is balanced by Prop. 1.3 in loc.cit. It follows that the homological equivalence coincides with the algebraic one for codimension $2$ algebraic cycles on $\bcC $, see \cite[Theorem 1(ii)]{BS}. The group $H^4(\bcC ,\ZZ )$ is isomorphic to $\ZZ $ by the Lefschetz's hyperplane section theorem and the Poincar\'e duality. This gives that
  $$
  CH^2(\bcC )=A^2(\bcC )\oplus \ZZ \; ,
  $$
and the action induced by $\tau $ on $CH^2(\bcC )$ splits into the action on $A^2(\bcC )$ and the identity action on $\ZZ $.

For any $i$ let $H^i(\bcC )$ be the cohomology of the complex cubic $\bcC $ with coefficients in $\QQ $. Recall that
  $$
  H^1(\bcC )=H^5(\bcC )=0\; ,\quad H^2(\bcC )=H^4(\bcC )
  =\QQ \quad \hbox{and}\quad H^3(\bcC )=\QQ ^{\oplus 10}\; .
  $$
As to Dolbeault cohomology, we have that $h^{3,0}(\bcC )=0$ and $h^{2,1}(\bcC )=5$.

Let $B_i$ be the pre-image of $\bcP _i$ in $\bcP $ under the above isogeny $\Lambda $. The Prymian $\bcP $ is generated by $B_2$ and $B_3$. Identifying $\bcP $ with $A^2(\bcC )$, and $A^2(\bcC )$ with the intermediate Jacobian $J^2(\bcC )$ via the Abel-Jacobi isomorphism, we can also look at $B_2$ and $B_3$ as two subgroups generating $A^2(\bcC )$ or $J^2(\bcC )$ respectively.

The genera of the curves $C_2'$ and $C_3'$ are $2$ and $4$ respectively, whence $\bcP $ is an abelian surface and $\bcP _3$ is an abelian threefold over $\CC $. Looking at the intermediate Jacobian of $\bcC $ as the quotient
  $$
  H^{2,1}(\bcC )^{\vee }/H_3(\bcC ,\ZZ )
  $$
and taking into account that any isogeny induces an isomorphism on the level of tangent spaces, we see that $\Lambda $ induces the splitting
  $$
  H^{2,1}(\bcC )=W_2\oplus W_3\; ,
  $$
such that the dual vector spaces $W_2^{\vee }$ and $W_3^{\vee }$ project from the tangent space to the intermedian Jacobian onto the groups $B_2$ and $B_3$ in $J^2(\bcC )$.

\begin{theorem}
\label{cubic action}
The involution $\tau ^*:A^2(\bcC )\to A^2(\bcC )$ acts identically on $B_3$, and as the multiplication by $-1$ on $B_2$. Similarly, the induced action on $H^{2,1}(\bcC )$ splits into the identity action on $W_3^{\vee }$ and multiplication by $-1$ on $W_2^{\vee }$.
\end{theorem}

\begin{pf}
The involution $\tau $ on the cubic $\bcC $ induces the involutions $\tau _2$ on $C'_2$ and $\tau _3$ on $C'_3$. In turn, they induce the involutions $\tau _2^*$ and $\tau _3^*$ on $J_2$ and $J_3$ respectively. Let $P$ be a point on the plane $\Pi _{\tau }$, let $\Pi _P$ be the span of $P$ and $l_{\tau }$, and look at the equation of the cubic $\bcC $,
  $$
  l_{00}(x_2,x_3,x_4)x_0^2+l_{11}(x_2,x_3,x_4)x_1^2+
  l_{01}(x_2,x_3,x_4)x_0x_1 +f_3(x_2,x_3,x_4)=0\; ,
  $$
see (\ref{S3+}) above. Under an appropriate change of the coordinates $x_3$ and $x_4$, keeping the coordinates $x_0$, $x_1$ and $x_2$ untouched, the plane $\Pi _P$ will be given by the equation
  $$
  \Pi _P : x_3=x_4=0\; .
  $$
Herewith, as the coordinates $x_0$, $x_1$ and $x_2$ remain the same, the involution $\tau $ in $\PR ^4$ can be expressed by the same formula, so that the equations for $l_{\tau }$ and $\Pi _{\tau }$ remain the same too (see Section \ref{idaction}). Substituting $x_3=x_4=0$ into the above equation for the cubic $\bcC $, we obtain the equation for the fibre
  $$
  \Pi _P\cap \bcC
  $$
of the projection
  $$
  p:\bcC \dashrightarrow \Pi _{\tau }
  $$

\medskip

\noindent over the point $P=(0:0:1:0:0)$ of the intersection of two planes $\Pi _P$ and $\Pi _{\tau }$. Namely, $\Pi _P\cap \bcC $ is given by the equation
  $$
  x_2(\alpha x_0^2+\beta x_1^2+\gamma x_0x_1+\delta x_2^2)=0\; ,
  $$
where $\alpha $, $\beta $, $\gamma $ and $\delta $ are some numbers in $\CC $. If a point $Q=(a_0:a_1:a_2:0:0)$ in $\Pi _P\cap \bcC $ is such that $a_2=0$ then $Q$ sits on the line $l_{\tau }$. As we are interested in the fibre of the projection from $\bcC \smallsetminus l_{\tau }$ we must set $x_2\neq 0$. Then, if $C_P$ is the Zariski closure of the set $(\Pi _P\cap \bcC  )\smallsetminus l_{\tau }$ in $\bcC $, the curve $C_P$ is the conic defined by the equation
  $$
  \alpha x_0^2+\beta x_1^2+\gamma x_0x_1+\delta x_2^2
  =0
  $$
in $\Pi _P$.

The point $P$ is in $C$ if and only if $C_P=l_P^++l_P^-$. Moreover,
  $$
  P\in C_3 \Leftrightarrow \delta =0\qquad \hbox{and}
  \qquad P\in C_2\smallsetminus C_3 \Leftrightarrow \delta \neq 0\; .
  $$
Then we see that, if $P\in C_3$, the lines $l_P^+$ and $l_P^-$ meet the plane $\Pi _{\tau }$ at the point $P$. It follows then that $\tau _3^*=\id $.

Suppose $P\in C'_2$. Since $C_P$ splits,
  $$
  \alpha x_0^2+\beta x_1^2+\gamma x_0x_1+\delta x_2^2=
  \delta (x_2+b_0x_0+b_1x_1)(x_2-b_0x_0-b_1x_1)\; ,
  $$
so that the lines of $C_P$ are defined by the equations
  $$
  l_P^+ : x_2+b_0x_0+b_1x_1=0 \qquad \hbox{and}
  \qquad l_P^- : x_2-b_0x_0-b_1x_1=0\; ,
  $$
which show that $\tau (l_P^+)=l_P^-$ and $\tau (l_P^-)=l_P^+$.

Thus, we obtain that the involution $\tau _2^*$ coincides with the involution $\iota _2^*$ on $J_2$, while the involution $\tau _3^*$ is the identity on $J_3$. This means that the action of $\tau ^*$ on the Prymian $\bcP _2$ is the multiplication by $-1$, and the action of $\tau ^*$ on the Prymian $\bcP _3$ is the identity. It follows that $\tau $ acts as multiplication by $-1$ on $B_2$ and identically on $B_3$. Looking at the tangent spaces, we also claim that $\tau $ induces the multiplication by $-1$ on $W_2^{\vee }$ and the identity on $W_3^{\vee }$.
\end{pf}

\medskip

Let now $\bcQ _1$ be yet another nonsingular $\tau $-invariant quadric in $\bcL _2$ and consider the pencil $\bcC \dasharrow \PR ^1$, which is the restriction of the linear system $|\bcQ _t|_{t\in \PR ^1}$ spanned by $\bcQ _0$ and $\bcQ _1$ onto $\bcC $. Let $\eta $ be the generic point of $\PR ^1$, $\bcC _{\eta }$ the generic fibre and $\bcC _{\bar \eta }$ the geometric generic fibre of the pencil $\bcC \dasharrow \PR ^1$. There exists a countable subset $Z$ in $\PR ^1$, such that the fibre $\bcC _{\bar \eta }$ is isomorphic, as a scheme over a subfield in $\CC (\PR^1)$, to the closed fibre $\bcC _t$, for all $t$ in $\PR ^1\smallsetminus Z$. Actually, $Z$ is the collection of all points with algebraic coordinate in $\AF ^1$ and the point at infinity $\infty $. The isomorphisms between $\bcC _{\bar \eta }$ and $\bcC _t$, for $t\in \PR ^1\smallsetminus Z$, commute with the action of the involution $\tau $. Then, by Theorem \ref{identity}, the action of $\tau ^*$ on $A^2(\bcC _{\bar \eta })$ is the identity.

Let $g:\bcC _{\eta }\to \bcC $ be the scheme-theoretical morphism of the generic fibre $\bcC _{\eta }$ into the cubic $\bcC $. Then $g$ induces the pull-back homomorphism $g^*:A^2(\bcC )\to A^2(\bcC _{\eta })$. This homomorphism is surjective, because any $\eta $-rational point on the surface $\bcC _{\eta }$ has its closure in the scheme $\bcC $.

Now, for any abelian group $A$ let $A_{\QQ }$ be the tensor product of $A$ with $\QQ $ over $\ZZ $. Theorem \ref{cubic action} implies that $B_2\cap B_3$ is a $2$-torsion subgroup in $A^2(\bcC )$. Therefore $A^2(\bcC )_{\QQ }$ is the direct sum of ${B_2}_{\QQ }$ and ${B_3}_{\QQ }$. By the same theorem, the involution $\tau ^*$ acts identically on ${B_3}_{\QQ }$ and as multiplication by $-1$ on ${B_2}_{\QQ }$.

Let $g^*_{\QQ }$ be the homomorphism induced by $g^*$ after tensoring with $\QQ $. The rational Chow group $A^2(\bcC _{\eta })_{\QQ }$ is embedded into $A^2(\bcC _{\bar \eta })$. The action of $\tau ^*$ on $A^2(\bcC _{\bar \eta })$ is the identity. This is why the action of $\tau ^*$ on $A^2(\bcC _{\eta })_{\QQ }$ is the identity too. Since $\tau ^*$ acts as multiplication by $-1$ on ${B_2}_{\QQ }$, we obtain that $g^*_{\QQ }({B_2}_{\QQ })=0$.

It means, in particular, that the rational Chow group $A^2(\bcC _{\eta })_{\QQ }$ of the $K3$-surface $\bcC _{\eta }$ over $\CC (\PR ^1)$ is covered by the $\tau ^*$-invariant component ${B_2}_{\QQ }$ of $A^2(\bcC )_{\QQ }$, which is isomorphic to the $\QQ $-localized Prymian ${\bcP _3}_{\QQ }$.


\bigskip

\medskip

\begin{small}

\end{small}

\bigskip

\bigskip

\begin{small}

{\sc Department of Mathematical Sciences, University of Liverpool, Peach Street, Liverpool L69 7ZL, England, UK}

\end{small}

\medskip

\begin{footnotesize}

{\it E-mail address}: {\tt vladimir.guletskii@liverpool.ac.uk}

\end{footnotesize}

\bigskip

\begin{small}

{\sc Department of Mathematics, Yaroslavl State University, 108 Respublikanskaya str., Yaroslavl 150000, RUSSIA}

\end{small}

\medskip

\begin{footnotesize}

{\it E-mail address}: {\tt astikhomirov@mail.ru, tikhomir@yspu.yar.ru}

\end{footnotesize}

\end{document}